\documentclass{amsart}
\usepackage{amssymb}
\usepackage{amsmath}
\usepackage{amscd}
\usepackage{amsbsy}
\usepackage{comment}
\usepackage{mathabx}
\usepackage{MnSymbol}
\usepackage{extarrows}
\usepackage{rotating}
\usepackage{mathrsfs}
\usepackage{lscape}
\usepackage[OT2, T1]{fontenc}
\usepackage[all,cmtip]{xy}

\DeclareMathOperator{\PGL}{PGL}
\DeclareMathOperator{\PSL}{PSL}

\DeclareMathOperator{\GL}{GL}
\DeclareMathOperator{\SL}{SL}

\DeclareMathOperator{\Ker}{Ker}

\DeclareMathOperator{\Gal}{Gal}

\newcommand{\Q}{{\mathbb Q}}

\newcommand{\C}{{\mathbb C}}

\newcommand{\F}{{\mathbb F}}
\newcommand{\PP}{{\mathbb P}}

\newcommand{\ff}{\mathfrak{f}}

\begin {document}

\newtheorem{thm}{Theorem}
\newtheorem{lem}{Lemma}[section]
\newtheorem{prop}[lem]{Proposition}

\newtheorem{cor}[lem]{Corollary}

\theoremstyle{definition}

\theoremstyle{remark}

\title[]{Base change for Elliptic Curves over Real Quadratic Fields}

\author{Luis Dieulefait and Nuno Freitas}
\address{
 Mathematisches Institut\\
Universit\"{a}t Bayreuth\\
95440 Bayreuth, Germany
}
\email{nunobfreitas@gmail.com}

\address{}

\email{ldieulefait@ub.edu}

\date{}
\thanks{The first-named author was supported by the MICINN grant MTM2012-33830 and ICREA Academia Research Prize. The second-named author was supported
through a grant within the framework of the DFG Priority Programme 1489
{\em Algorithmic and Experimental Methods in Algebra, Geometry and Number Theory}.
}


\begin{abstract} Let $E$ be an elliptic curve over a real quadratic field $K$ and $F/K$
a totally real finite Galois extension. We prove that $E/F$ is modular.
\end{abstract}
\maketitle

\makeatletter
\renewcommand{\tocsection}[3]{%
  \indentlabel{\@ifnotempty{#2}{\ignorespaces#1 #2.\quad}}#3\dotfill
}

\makeatother

\section{Introduction}

For $F$ a totally real number field we write $G_F := \Gal(\overline{\Q}/F)$ for its absolute Galois group.
For a Hilbert modular form $\ff$ we denote by $\rho_{f,\lambda}$ its attached $\lambda$-adic representation.
We say that a continuous Galois representation $\rho : G_F \rightarrow \GL_2(\overline{\Q}_\ell)$ is {\bf modular}
if there exists a Hilbert newform $\ff$ and a prime $\lambda \mid \ell$ in its field of coefficients
$\Q_\ff$ such that we have an isomorphism $\rho \sim \rho_{\ff, \lambda}$.
In \cite{Di2} and \cite[Section 5]{Di}, the first named author proved base change for the $\GL_2$ case over $\Q$ (\cite[Theorem 1.2]{Di}).

\begin{thm} Let $f$ be a classical modular form of weight $k \geq 2$ and field of coefficients $\Q_f$. For a prime $\lambda$ of $\Q_f$ write $\rho_{f,\lambda}$ for the attached $\lambda$-adic representation. Let $F/\Q$ be a totally real number field. Then the Galois representation $\rho_{f,\lambda}|G_F$ is (Hilbert) modular in the sense above.
\label{thm:basechange}
\end{thm}

In the recent paper \cite{FHS} the following modularity theorem is proved.
\begin{thm} \label{thm:modularity}
Let $E$ be an elliptic curve defined over a real quadratic field $K$.
Then $E$ is Hilbert modular over $K$. 
\end{thm}

The aim of this note is to establish a base change result for certain elliptic curves
as a consequence of Theorem~\ref{thm:modularity}. More precisely, we prove the following.

\begin{thm} Let $E$ be an elliptic curve over a real quadratic field $K$. Let also $F/K$ be
a totally real finite Galois extension. Then $E/F$ is modular.
\label{thm:main}
\end{thm}

This result has applications in the context of the Birch and Swinnerton-Dyer conjecture. 
Indeed, modularity of $E$ after base change guarantees that the $L$-function $L(E/F,s)$
is holomorphic in $\C$ and, in particular, its order of vanishing at $s=1$ is a well
defined non-negative integer, in agreement with what is predicted by the BSD conjecture.
Furthermore, modularity of $E/F$ allows the construction of Stark-Heegner points on $E$ over (not necessarily real) quadratic extensions of $F$. For details regarding this application we refer the reader to \cite{GRZ} and the references therein.

\bigskip

\noindent \textbf{Acknowledgements.} We would like to thank K\k{e}stutis \v{C}esnavi\v{c}ius, Jos{\'e} Mar{\'i}a Giral, Victor Rotger and Samir Siksek for their useful comments. We also thank the anonymous referee for his comments.
 
\section{Elliptic curves with big non-solvable image mod $p=3,5$ or $7$}

Let $F / K$ be a finite extension of totally real number fields. Let $E/K$ be an elliptic curve. We will say that $\overline{\rho}_{E,p}(G_F)$ {\bf is big} if $\overline{\rho}_{E,p}(G_{F(\zeta_p)})$ is absolutely irreducible, otherwise we say it {\bf is small}. In particular, if $\overline{\rho}_{E,p}(G_{F})$ is non-solvable then it is big. We now restate \cite[Theorem 3 and 4]{FHS}.

\begin{thm}\label{thm:357}
Let $p=3,5$ or $7$. Let $F / K$ and $E/K$ be as above.
Suppose that $\overline{\rho}_{E,p}(G_F)$ is big.
Then $E$ is modular over $F$.
\end{thm}

The following proposition is well-known.

\begin{prop} Let $F / K$ be a finite Galois extension of totally real fields and $E/K$ an elliptic curve. Let $p$ be a prime and suppose that $\overline{\rho}_{E,p}(G_K)$ is non-solvable. Then $\overline{\rho}_{E,p}(G_F)$ is non-solvable.
\label{pp:nsolv}
\end{prop}
\begin{proof} Since $\overline{\rho}_{E,p}(G_K)$ is non-solvable we have $p > 3$. From Dickson's theorem (see also Proposition~\ref{prop:dickson}), having $\overline{\rho}_{E,p}(G_K)$ non-solvable implies that projectively $\overline{\rho}_{E,p}(G_K)$ is $A_5$ or $\PSL_2(\F_p)$ or $\PGL_2(\F_p)$. For the last
two cases the proposition is a particular case of \cite[Lemma 3.2]{Di2}. Since $A_5$ is a simple group
the same argument as in \textit{loc. cit.} also applies in this case.
\end{proof}

We have the following corollary.

\begin{cor}\label{cor:nonsolv} Let $F / K$ and $E/K$ be as in Proposition~\ref{pp:nsolv}. Let $p=3,5$ or $7$.
Suppose that $\overline{\rho}_{E,p}(G_K)$ is non-solvable. Then $E$ is modular over $F$.
\end{cor}

\begin{proof} From the previous proposition we have that $\overline{\rho}_{E,p}(G_F)$ is non-solvable, hence it is big. Thus $E/F$ is modular by Theorem~\ref{thm:357}.
\end{proof}

\section{Elliptic curves with projective image $S_4$ or $A_4$ mod $p=3,5$ or $7$}

Let $E/K$ be an elliptic curve. We have seen that if $\overline{\rho}_{E,p}$ has big non-solvable image then after a base change to a totally real Galois extension its image is still non-solvable. We now want to understand what can happen when $\overline{\rho}_{E,p}(G_K)$ is big and solvable. We first recall the following well know fact.

\begin{prop} Let $E/K$ be an elliptic curve. Write $G$ for the image of $\overline{\rho}_{E,p}$ in $\GL_2(\F_p)$ and $H$ for its image in $\PGL_2(\F_p)$. Then, there are the following possibilities:
\begin{enumerate}
\item[(a)] $G$ is contained in a Borel subgroup;
\item[(b)] $G$ contains $SL_2(\F_p)$;
\item[(c)] $H$ is cyclic, $G$ is contained in a Cartan subgroup;
\item[(d)] $H$ is dihedral, $G$ is contained in the normalizer of a Cartan subgroup;
\item[(e)] $H$ is isomorphic to $A_4$, $S_4$ or $A_5$
\end{enumerate}
\label{prop:dickson}
\end{prop}

Let $p=3,5$ or $7$. Let also $G$ and $H$ be as in the proposition. Remembering that
$\PSL_2(\F_\ell)$ is simple for $p \geq 5$, by Jordan-Moore's theorem, and that $\PSL_2(\F_3) \simeq S_4$, we divide the cases where $\overline{\rho}_{E,p}(G_K)$ is big and solvable into two types:
\begin{enumerate}
 \item[(I)]  $H \cong S_4$ or $A_4$,
 \item[(II)] $H$ is dihedral.
\end{enumerate}
Suppose we are in case (I). Let $F / K$ be a finite Galois extension and set $H_F :=\PP(\overline{\rho}_{E,p}(G_F))$. We would like that $H_F$ is also isomorphic to $A_4$ or $S_4$ since this would mean that $\overline{\rho}_{E,p}(G_F)$ is big and Theorem~\ref{thm:357} applies. Since $F / K$ is Galois we have that $H_F$ is normal subgroup of $H$. Write $I = \{1 \}$ for the trivial group and $D_4$ for the dihedral group in four elements. The normal subgroups of $S_4$ and $A_4$ are respectively
\begin{itemize}
 \item $I$, $D_4$, $A_4$ and $S_4$,
 \item $I$, $D_4$ and $A_4$.
\end{itemize}
Thus, the cases where Theorem~\ref{thm:357} does not apply over $F$ are when the pair of groups $(H, H_F)$ is one of
\begin{equation}\label{eq:cases}
(S_4, D_4), \quad (S_4, I), \quad  (A_4, D_4), \quad (A_4, I).
\end{equation}
Since we are working with totally real fields the complex conjugation has projective image of order 2. Thus the cases with $H_F = I$ cannot happen.

\subsection{A Sylow base change} We now deal with the remaining cases from \eqref{eq:cases}. Recall that we want to base change $E/K$ to $F$ where $F/K$ is finite and Galois. Suppose that $(H, H_F)$ is $(S_4, D_4)$ or $(A_4, D_4)$. Let $F_3$ be a subfield of $F$ such that the Galois group $\Gal(F/F_3)$ is a $3$-Sylow subgroup of $\Gal(F/K)$. In particular, $F/F_3$ is a solvable extension. We shall shortly prove the following.

\begin{lem}\label{lem:sylow} The projective image $H_{F_3} := \PP(\overline{\rho}_{E,p}(G_{F_3}))$ is isomorphic to $S_4$ or $A_4$. In particular, $\overline{\rho}_{E,p}(G_{F_3})$ is big.
\end{lem}

\noindent From this lemma and Theorem~\ref{thm:357} it follows that $E/F_3$ is modular. Finally, an application of Langlands solvable base change (see \cite{Langlands}) allows to conclude that $E/F$ is modular.

\bigskip

For the proof of Lemma~\ref{lem:sylow} we will need the following
elementary lemma from group theory.

\begin{lem} Let $G$ be a profinite group. Let $M \subset G$ be a subgroup of finite index $i$. Let $N$ be a normal subgroup of $G$.
Write $j$ for the index of $M/(N \cap M)$ in $G/N$. Then $ j \mid i$.
\label{lem:index}
\end{lem}

\begin{proof}
We prove it for the case of finite groups. The required divisibility
follows from the following elementary equalities:
$$ |G|  =  |N| \cdot  [G:N] $$
$$  |M|  =  |N \cap M | \cdot [M: N \cap M] $$
Dividing the first equality into the second, we conclude that $j$
divides $i$.
\end{proof}

\begin{proof}[Proof of Lemma~\ref{lem:sylow}] Let $F_3$ be as above and set
\[
G :=\Gal(\overline{\Q}/K), \quad M :=\Gal(\overline{\Q}/F_3), \quad N:=\Ker(\PP\overline{\rho}_{E,p}).
\]
Let $L/K$ be the Galois extension fixed by $N$. Observe that $L/L\cap F_3$ is Galois and
\[
 G/N \cong \Gal(L/K), \qquad M/(M\cap N) \cong \Gal(L/L\cap F_3).
\]
From Lemma~\ref{lem:index} we see that
\[
[\Gal(L/K) : \Gal(L/L\cap F_3)] = j \mid i = [G : M]
\]
and we also have
\[
 |\Gal(L/K)| = j|\Gal(L/L\cap F_3)|.
\]
Note that $\Gal(L/L\cap F_3) \cong H_{F_3}$. From the way we choose
$F_3$ it is clear that $3 \nmid i$, hence $3 \nmid j$. By
hypothesis $G/N \cong S_4$ or $A_4$, hence 3 divides $|\Gal(L/K)|$
and $|H_{F_3}|$. Finally, the conditions $3 \mid |H_{F_3}|$ and $D_4
\subset H_{F_3}$ together imply that $H_{F_3}$ is isomorphic to
$S_4$ or $A_4$.
\end{proof}

We summarize this section into the following corollary.

\begin{cor}\label{cor:solv} Let $F/K$ be a finite Galois extension of totally real fields. Let $E/K$ be an elliptic curve.
Suppose that for $p=3,5$ or $7$ we have that $\overline{\rho}_{E,p}(G_K)$ is big and solvable. Suppose further that
$\PP(\overline{\rho}_{E,p}(G_K)) \cong S_4$ or $A_4$. Then $E/F$ is modular.
\end{cor}

Everything we have done so far works for any Galois extension $F/K$. Moreover, it is clear that the remaining
cases are those when $\overline{\rho}_{E,p}(G_K)$ is small or projectively dihedral simultaneously for $p=3,5,7$. The restriction in
the statement of Theorem~\ref{thm:main} to quadratic fields arises precisely from dealing with them, which is
the content of the next section.

\section{Elliptic curves having small or projective Dihedral image at $p=3,5$ and $7$}
\label{sec:smallimage}

Let $K$ be a real quadratic field. From Theorem~\ref{thm:357} an elliptic curve $E/K$ is modular over $K$ except possibly if $\overline{\rho}_{E,p}(G_K)$ is small simultaneously for $p=3$, $5$, $7$. Suppose $K \neq \Q(\sqrt{5})$. In \cite{FHS} it is shown that such an elliptic curve gives rise to a $K$-point on one of the following modular curves:
\begin{gather*}
X(\mathrm{b}5,\mathrm{b}7), \qquad
X(\mathrm{b}3,\mathrm{s}5), \qquad
X(\mathrm{s}3,\mathrm{s}5),\\
\label{eqn:last4}
X(\mathrm{b}3,\mathrm{b}5,\mathrm{d}7), \quad
X(\mathrm{s}3,\mathrm{b}5,\mathrm{d}7), \quad
X(\mathrm{b}3,\mathrm{b}5,\mathrm{e}7), \quad
X(\mathrm{s}3,\mathrm{b}5,\mathrm{e}7),
\end{gather*}
where $\mathrm{b}$ and $\mathrm{s}$
respectively stand for \lq Borel\rq and \lq normalizer of split Cartan\rq. The notation $\mathrm{d}7$ and $\mathrm{e}7$ is explained in \cite[section 10]{FHS}, here we remark only that they indicate mod $7$ level structures that are respectively finer than \lq normalizer of split Cartan\rq\
and \lq normalizer of non-split Cartan\rq.
Denote by $\mathcal{E}_K$ the set of elliptic curves (up to quadratic twist) corresponding to $K$-points in the previous modular curves.

\smallskip

\noindent In \cite{FHS} it is also shown that an elliptic curve $E/\Q(\sqrt{5})$ with simultaneously small image for $p=3,5,7$ gives rise to a $\Q(\sqrt{5})$-point in one of the following modular curves
\begin{equation*}
X(\mathrm{d}7), \quad X(\mathrm{e}7), \quad X(\mathrm{b}3,\mathrm{b7}), \quad X(\mathrm{s}3,\mathrm{b}7).
\end{equation*}
Denote by $\mathcal{E}_{\Q(\sqrt{5})}$ the set of elliptic curves (up to quadratic twist) corresponding to $\Q(\sqrt{5})$-points in these four modular curves.

\bigskip
Furthermore, it is also follows from \cite{FHS} that, for any real quadratic field $K$, we have
\begin{itemize}
 \item[(i)] $\mathcal{E}_K$ contains all elliptic curves (up to quadratic twist) with small or projective Dihedral image simultaneously at $p=3,5,7$;
 \item[(ii)] $\mathcal{E}_K$ is finite;
 \item[(iii)] Let $E \in \mathcal{E}_K$. Then, either $E$ is a $\Q$-curve or $E$ has complex multiplication or $\bar{\rho}_{E,7}(G_K)$ contains $\SL_2(\F_7)$.
\end{itemize}

We can now easily prove the following.
\begin{cor}\label{cor:simage} Let $K$ be a real quadratic field. Let $E \in \mathcal{E}_K$. Let $F/K$ be a finite totally real Galois extension. Then $E / F$ is modular.
\end{cor}
\begin{proof} From (iii) above we know that either (a) $E/K$ is a $\Q$-curve or has complex multiplication or (b) $\rho_7(G_K)$ is non-solvable. Suppose we are in case (a). Base change follows from \cite[Proposition 12.1]{JacLang} in the CM case; if $E$ is a $\Q$-curve, by results of Ribet and Serres's conjecture (now a theorem due to Khare-Wintenberger), it arises from a classical modular form thus base change follows by Theorem~\ref{thm:basechange}. In case (b), it follows from Corollary~\ref{cor:nonsolv} that $E/F$ is modular.
\end{proof}

\section{Proof of the main theorem}

Let $K$ be a real quadratic field and $E/K$ an elliptic curve.  Write $\bar{\rho}_p = \bar{\rho}_{E,p}$. The curve $E/K$ must satisfy at least one of the following three cases
\begin{enumerate}
 \item $\bar{\rho}_p(G_K)$ is big and non-solvable for some $p \in \{3,5,7\}$,
 \item $\bar{\rho}_p(G_K)$ is big, solvable and satisfy $\PP(\bar{\rho}_p(G_K)) \cong S_4, A_4$ for some $p \in \{3,5,7\}$,
 \item $E/K$ belongs to the set $\mathcal{E}_K$.
\end{enumerate}
\medskip

\noindent Let $F/K$ be a totally real finite Galois extension. In each case, modularity of $E/F$ now follows directly from one
of the previous sections:

Case (1): this is Corollary~\ref{cor:nonsolv}.

Case (2): this is Corollary~\ref{cor:solv}.

Case (3): this is Corollary~\ref{cor:simage}.
\qed

\end{document}